\documentclass[12pt]{amsart}
\usepackage{amsmath,amssymb, xypic,
verbatim,amscd}
\numberwithin{equation}{section}

\newcommand\pone{\Bbb{P}^1}
\newcommand\mm{{m}}

\newcommand\mf{\mathcal{F}}
\newcommand\spanA{\operatorname{span}_A}
\newcommand\ord{\operatorname{ord}}
\newcommand\rk{\operatorname{rk}}

\newcommand\tiTheta{\tilde{\Theta}}

\newcommand{\spec}{\operatorname{Spec}}
\newcommand{\home}{\operatorname{Hom}}
\newcommand\pt{\operatorname{pt}}
\newcommand\aA{\text{\bf {a}}}

\newcommand\red{\operatorname{red}}
\newcommand\mv{\mathcal{V}}
\newcommand\cee{\Bbb{C}}
\newcommand\zee{\Bbb{Z}}
\newcommand\C{\Bbb{C}}
\newcommand\ma{\mathcal{A}}

\newcommand\tensor{\otimes}
\newcommand\Wr{\operatorname{Wr}}
\newcommand\mw{\mathcal{W}}

\newcommand\Fl{\operatorname{Fl}}

\newcommand\Gr{\operatorname{Gr}}

\newcommand\mq{\mathcal{Q}}
\newcommand\mi{\mathcal{I}}

\newcommand\sll{\operatorname{sl}}
\newcommand\emm{\mathfrak{m}}
\newcommand\pp{\mathfrak{p}}
\newcommand\spanC{\operatorname{span}_{\cee}}

\newcommand\bull{\sssize{\bullet}}

\newtheorem{theorem}{Theorem}[section]
\newtheorem{lemma}[theorem]{Lemma}

\newtheorem{proposition}[theorem]{Proposition}
\newtheorem{corollary}[theorem]{Corollary}
\newtheorem{definition-proposition}[theorem]{Definition-Proposition}
\theoremstyle{definition} 
\theoremstyle{defi}

\theoremstyle{remark} \newtheorem{remark}[theorem]{Remark}

\newcommand{\bs}{\boldsymbol}

\begin{document}
\title[Critical Points and Schubert calculus]
{Multiplicity of critical points of master functions and Schubert calculus}
\author[{}]{ Prakash Belkale${}^{*\ 1}$, 
Evgeny Mukhin ${}^{**\ 2}$,
\and Alexander Varchenko${}^{***\ 3}$}

\thanks{${}^1$ Supported in part by NSF grant  DMS-0300356}

\thanks{${}^2$ Supported in part by NSF grant DMS-0140460}

\thanks{${}^3$ Supported in part by NSF grant DMS-0244579}

\begin{abstract}
In ~\cite{mukvar}, some correspondences were defined between critical points
of master functions associated to $\sll_{N+1}$ and subspaces of
$\cee[x]$ with given ramification properties. In this paper we show
that these correspondences are fact scheme theoretic isomorphisms of
appropriate schemes. This gives relations between multiplicities of
critical point loci of the relevant master functions and
multiplicities in Schubert calculus.
\end{abstract}
\maketitle
\centerline{\it ${}^{*\ ***}$Department of Mathematics, University of
  North Carolina at Chapel Hill,} \centerline{\it Chapel Hill, NC
  27599-3250, USA} \medskip 
\medskip \centerline{\it ${}^{**}$
Department of Mathematical Sciences,} 
\centerline{\it Indiana University 
Purdue University Indianapolis,}
\centerline{\it 402 North Blackford St., Indianapolis,
IN 46202-3216, USA}
 \medskip

\section*[intro]{Introduction}
In ~\cite{mukvar}, a correspondence between the  following objects was shown:
\begin{enumerate}
\item 
Critical points of certain master functions associated to the Lie
algebra $\sll_{N+1}$;
\item Vector subspaces $V\subset \cee[x]$ of 
rank $N+1$ with prescribed ramification at points of $\cee\cup\{\infty\}$.
\end{enumerate}
Let $l_1, \dots , l_N$ be nonnegative integers and $z_1,\dots,z_n$
distinct complex numbers.  For $s=1, \dots , n$, fix nonnegative integers
$m_s(1), \dots , m_s(N)$. Define polynomials 
$T_1, \dots , T_N$ by the formula
$T_i\,  =\,  \prod_{s=1}^n \, (x-z_s)^{m_{s}(i)}$. 
{\it The master function $\Phi$ associated to this data} is the rational function 
$$
\Phi(t^{(i)}_j)\ =\ \prod_{i=1}^N \prod_{j=1}^{l_i} \ T_i(t^{(i)}_j)^{-1} \
\prod_{i=1}^{N-1} \prod_{j=1}^{l_i}\prod_{k=1}^{l_{i+1}}\
(t^{(i)}_j- t^{(i+1)}_k)^{-1}\
\prod_{i=1}^{N}\prod_{1\leq j< k \leq l_i}\
(t^{(i)}_j-t^{(i)}_k)^{2}
$$ 
of $l_1+\dots+l_N$ variables $t^{(1)}_1,\dots,t^{(1)}_{l_1},\dots,
t^{(N)}_{1},\dots,t^{(N)}_{l_N}$ considered on the set of points where
\begin{enumerate}
\item[$\bullet$] The numbers $t^{(i)}_1, \dots , t^{(i)}_{l_i}$ are distinct;
\item[$\bullet$]
The sets $\{t^{(i)}_1, \dots , t^{(i)}_{l_i}\}$
and $\{t^{(i+1)}_1, \dots , t^{(i+1)}_{l_{i+1}}\}$ do not intersect;
\item[$\bullet$]
The sets $\{t^{(i)}_1, \dots , t^{(i)}_{l_i}\}$
and $\{z_1, \dots , z_n\}$ do not intersect.
\end{enumerate}
The master functions considered in (1) are functions $\Phi$ as
above. These functions appear in hypergeometric solutions to the KZ
equation. They also appear in the Bethe ansatz method of the Gaudin model,
where the goal is to write formulas for singular vectors in a
tensor product of representations of $\sll_{N+1}$ 
 starting from a critical point. Those Bethe vectors are
eigenvectors of certain commuting linear operators called Hamiltonians 
of the Gaudin model. We refer the reader
to ~\cite{var} for a detailed discussion of these themes.

The set of objects given by (2) are important in the study of linear
series on compact Riemann surfaces (see ~\cite{eisenharris}).

In both of the objects (1), (2) above, there is a natural notion of
multiplicity. In (1), we can consider the geometric multiplicity of
the critical scheme (see Proposition
~\ref{transversalmultiplicity}). In (2), we may view the set of such
objects as the intersection of some Schubert cells in a Grassmannian
$\Gr(N+1,\cee_{d}[x])$ and hence, there is an associated intersection
multiplicity (which is in this case the same as the geometric
multiplicity at the given $V$ of the intersection of Schubert
cells). In this paper, we show that these multiplicities agree
(Theorems ~\ref{mainte}, ~\ref{mainte2}, and corollaries ~\ref{MLK},
~\ref{ba11}).

One consequence of such an agreement, is that intersection numbers in
Grassmannians can be calculated from the critical scheme of master
functions and vice versa. As indicated in ~\cite{mukvar}, this relation provides a link via critical schemes of master functions, between representation theory of $\sll_{N+1}$ and Schubert calculus (see ~\cite{be} for a different geometric approach).

A related consequence is that the intersection of associated Schubert
cells is transverse at $V$ if and only if the associated critical
scheme is of geometric multiplicity $1$.

A note on our methods: We show the equalities of multiplicities by
using Grothendieck's scheme theory. To obtain statement on
multiplicities, we need to show that the correspondences are in fact
isomorphisms of (appropriate) schemes. By Grothendieck's functorial
approach to schemes, this aim will be achieved if we can replace
$\cee$ in ~\cite{mukvar} by an arbitrary local ring $A$ and show that
the correspondences hold with objects over $A$.

This requires us to develop the theory of Wronskian equations over an
arbitrary local ring, in particular to develop criteria for
solvability in a purely algebraic manner (see Lemma ~\ref{BA}). We
also need to revisit key arguments in ~\cite{mukvar} and modify their
proof so that they apply over any local ring (see Theorem
~\ref{Main}).

\subsection{Notation}
For a ring $A$, let $A_{ d}[x]$ denote the set of polynomials with
coefficients in $A$ of degree $\leq d$. An element in $A[x]$ is said
to be monic, if its leading coefficient is invertible. The
multiplicative group of units in $A$ is denoted by $A^*$. 

 For a ring $A$ and elements $y_1,\dots,y_k\in A$, we will denote by $(y_1,\dots,y_k)$, the ideal generated by the $y_i$ in $A$.

A local ring $(A,\emm)$ over $\cee$ is a Noetherian ring $A$ containing
$\cee$ with a unique maximal ideal $\emm$. The residue field of the
local ring is defined to be ${A/\emm}$. We will only consider local rings
containing $\cee$ with residue field $\cee$. A scheme in the paper
stands for an algebraic scheme over $\cee$.

We denote the permutation group of the set  $\{1,\dots,N\}$  by $\Sigma_{N}$.
\section{Formulation of the main result}
\subsection{Some preliminaries}\label{1.1}

For a vector bundle $\mw$ on a scheme $X$, denote by $\Fl(\mw)\to X$
the fiber bundle  whose fiber over a point $x\in X$ is the flag variety of
complete filtrations of the fiber $\mw_x$.

Let $W$ be a vector space of rank $d+1$ and $N$, $0\leq N \leq d$, an
integer. There is a natural exact sequence of vector bundles on the
Grassmannian $\Gr(N+1,W)$,
$$
0\ \to\ \mv \
\to\ W\tensor \mathcal{O} \to\ \mq\ \to\ 0\ .
$$
The fiber of this sequence at a point $V\in \Gr(N+1,W)$ is 
$$
0\ \to \ V\ \to\ W\ \to\ W/V\ \to\ 0\ .
$$ 

Let $\mf$ be a complete flag on $W$:
$$
\mf\ :\ \text{ }\{0\} = F_0
\ \subsetneq F_1\ \subsetneq \ \dots\ \subsetneq \ F_{d+1}\ =\ W\ .
$$ 
A {\it ramification sequence} is a sequence  $\aA$ of the form $(a_1,\dots,a_k)\in \zee^k$ such that $a_1 \geq \dots \geq a_{k} \geq 0.$ For a ramification sequence $\aA=(a_1,a_2,\dots,a_{N+1})$
satisfying ${a}_1\leq d-N$,
define the Schubert cell 
\begin{align}
&
\Omega^o_{\aA}(\mf) \ =\ \{V\in \Gr(N+1,W)\mid \rk(V\cap F_u) = \ell ,\
\notag
\\
&
\phantom{aaaaaaaaaa}
d-N+\ell-{a}_{\ell} \leq u < d - N + \ell + 1 - {a}_{\ell+1},
\ \text{}\ell = 0, \dots,N+1 \}\ ,
\notag
\end{align}
where ${a}_0 = d-N$,\ ${a}_{N+2} = 0$. The cell
$\Omega^o_{\aA}(\mf)$ is a smooth connected variety. 
The closure of $\Omega^o_{\aA}(\mf)$ is denoted by 
$\Omega_{\aA}(\mf)$. The codimension 
of $\Omega^o_{\aA}(\mf)$ is
$$
|\aA|\ =\ {a}_1\ +\ 
{a}_2\ +\ \dots\ +\ {a}_{N+1}\ . 
$$ 
Every 
$N+1$-dimensional vector subspace of $W$  belongs to a unique Schubert cell
$\Omega^o_{\aA}(\mf)$.

Denote by $\mv_{\aA}$ the pull back of $\mv$ to
$\Omega^o_{\aA}(\mf)\hookrightarrow\Gr(N+1,W)$.  There is
{\it the distinguished section} of 
$$
\Fl(\mv_{\aA})\ \to\
\Omega^o_{\aA}(\mf)
$$ 
which assigns to a point 
$V\in \Omega^o_{\aA}(\mf)$ the complete filtration
$$
0\ \subsetneq F_{d-N+1-a_1}\cap V\ 
\subsetneq F_{d-N+2-a_2}\cap V\ \subsetneq 
\ \dots \ \subsetneq\ F_{d-N+N+1-a_{N+1}}\cap V = V\ .
$$
This section may be used to partition each fiber  
$\Fl((\mv_{\aA})_{V})$ into Schubert cells (see definitions in Section ~\ref{fri}). 
This partition varies algebraically with $V$. 
That is, there is  a decomposition of $\Fl(\mv_{\aA})$ into relative Schubert cells 
each of which is a locally trivial (in the Zariski topology) 
fiber bundle over $\Omega^o_{\aA}(\mf)$.

\subsection{Schubert cells in flag varieties}\label{fri}
Let $V$ be a vector space of rank $N+1$, $\mf$ a complete flag on $V$
and $w\in\Sigma_{N+1}$. Define the Schubert cell $G^o_w(\mf)$
correponding to $w$ to be the subset of $\Fl(V)$ formed by points\
$ V_1\,\subsetneq V_2\,\subsetneq\,\dots\,\subsetneq \,V_{N+1}=V$\ such that
there exists a basis $u_1,\dots,u_{N+1}$ of $V$ satisfying the conditions
$$
V_i\ =\ \spanC\,(u_1,\dots u_i),\ \ u_i\in F_{w(i)},\ i=1,\dots, N+1\ .
$$

It is easy to see that 
$$
\Fl(V)\ =\ \bigsqcup_{w\in
\Sigma_{N+1}}G^o_w(\mf)\ .
$$ 
It is easy to see that the permutation of highest length
gives the {\em open cell} in this partition of $\Fl(V)$.


\subsection{Intersection theory in the space of functions}\label{secion}
Let $W=\cee_{ d}[x]$ be the space of polynomials of degree not
greater than $d$. Each point $z\in\cee \cup {\infty}$ determines a
full flag in $W$:
$$
\mf(z)\ : 0=F_0(z)\ \subsetneq \ F_1(z)\ 
\subsetneq\ \dots\ \subsetneq\  F_{d+1}(z) = W
$$
where for any $z\in \cee$ and $i$, $F_i(z)=(x-z)^{d+1-i}\cee[x]\cap W$,
and if $z=\infty$, then $F_i(z)$ is the  space of polynomials of degree $< i$.

For $V\in \Gr(N+1,W)$ and $z\in \cee\cup\infty$, there exists a unique
ramification sequence $\aA(z) = (a_1,\dots,a_{N+1})$, with 
$a_1\leq d-N$,
such that $V\in \Omega^o_{\aA(z)}(\mf(z))$. The sequence is called {\it the ramification
sequence of $V$ at $z$}. 

If $z\in\cee$, then this means that $V$ has a basis of the form
$$
(x-z)^{N+1-1+a_1}f_1,\ (x-z)^{N+1-2+a_2}f_2 ,\ \dots ,\
 (x-z)^{N+1-(N+1)+a_{N+1}}f_{N+1}
$$ 
with $f_i(z) \neq 0$ for $i=1,\dots,N+1$. The numbers $\{N+1-i+a_i\ |
\ i = 1, \dots , N+1\}$, are called the {\it exponents} of $V$ at $z$. 

If $z=\infty$, then the condition is that $V$ has a basis of the form $f_1,
f_2 , \dots , f_{N+1}$ with $\deg\,f_i\, = d-(N+1)+i-a_i$ for $i = 1,
\dots , N+1$. The numbers $\{d-(N+1)+i-a_i ,\ |\ i = 1, \dots , N+1\}$, are 
called the {\it exponents} of $V$ at $\infty$. 
\begin{remark}
The set of exponents of $V$ at any point in $\cee\cup\ \{\infty\}$ 
is a subset of $\{0, \dots , d\}$ of cardinality $N+1$. 
\end{remark}   
A point $z\in\cee\cup\{\infty\}$ is called {\it a ramification point of $V$}
if $\aA(z)$ is a sequence with at least one nonzero term.

For a finite dimensional subspace 
$E \subseteq \cee[x]$,  define the Wronskian $\Wr(E)
\in\cee[x]/ {\cee}^{*}$ as the Wronskian of a basis of $E$. 
The Wronskian of a subspace is a nonzero polynomial with the following
properties.

\begin{lemma}
 If\  $V\subseteq W$ lies 
in the   Schubert cell 
$\Omega^o_{\aA}(\mf(z))\subseteq \Gr(N+1,W)$ 
for some $z\in \cee$, then $\Wr(\, V)$  has a root at $z$ of multiplicity
$|\aA|$.


If \ $V\in  
\Omega^o_{\aA}(\mf(\infty))$, then $\deg\,\Wr(\, V)\, = (N+1)(d-N)-|\aA|$.
\end{lemma} 

We will fix the following objects:
\begin{enumerate}\label{sit}
\item[$\bullet$] A space of polynomials $W=\cee_{ d}[x]$;
\item[$\bullet$]
 A Grassmannian $\Gr(N+1,W)$ with the universal subbundle $\mv$;
\item[$\bullet$] 
Distinct points $z_1,\dots,z_n$ on $\Bbb{C}$;
\item[$\bullet$]
 At each point $z_s$, a ramification sequence $\aA(z_s)$ 
and at $\infty$, a ramification sequence $\aA(\infty)$
so that
\begin{equation}\label{eqn1}
\sum_{z\in\{z_1,\dots,z_n\}}|\aA(z)| +
|\aA(\infty)| =\dim\,\Gr(N+1,W)\,=\,(N+1)(d-N).
\end{equation}
\end{enumerate}
We set
\begin{align}
\Omega \ & =\ \bigcap_{s=1}^n\ \Omega^o_{\aA(z_s)}(\mf(z_s))\ \cap\
\Omega^o_{a{(\infty})}(\mf(\infty))\ , 
\notag
\\
K_i & =\ \prod_{s=1}^n \ (x-z_s)^{\sum_{\ell=1}^i a_{N+1-i+\ell}(z_s)}\ , 
\qquad  i = 1 , \dots , N+1\ ,  
\notag
\\
T_i &=\ K_{i+1} K_{i-1}/K_i^2 \ =\  
\prod_{s=1}^n \ (x-z_s)^{a_{N+1-i}(z_s)-a_{N+2-i}(z_s)}\ ,
\qquad  i = 1 , \dots , N\ .
\notag
\end{align}
We set $K_0 = 1$. 
Notice that $T_i$ is a polynomial.
  
The collection of these objects will be called the ``{\it basic
situation}''.

\begin{remark}
The following are basic facts from intersection theory in the space of polynomials.
\begin{itemize}
\item $\Omega$ is a finite scheme, its support is a 
finite set.  
\item In the definition of $\Omega$, 
we intersected Schubert cells. 
We obtain the same intersection if we intersect the closures of 
the same Schubert cells:
$$
\Omega\ =\
\bigcap_{s=1}^n \
 \Omega_{\aA(z_s)}(\mf(z_s))\ 
\cap \ \Omega_{\aA{(\infty)}}(\mf(\infty))\ .
$$
\end{itemize}
\end{remark}
Define $\Fl$ as a pull back 
\begin{equation}\label{cartsquare}
\xymatrix{\Fl\ar[r]\ar[d]^{\pi}  &                \Fl(\mv)\ar[d]\\
            \Omega\ar[r]                          & \Gr(N+1,W) }
\end{equation}
Points of $\Fl$ are pairs $(V,E_{\bull})$ where $V\in\Omega$ and 
$E_{\bull}=
(E_1\subsetneq E_2\subsetneq\dots\subsetneq E_{N+1}=V)$ is a 
complete filtration of $V$. For $V\in \Omega$, 
there are $n+1$ {\it distinguished
complete filtrations on} $V$ corresponding respectively
to the flags  $\mf(z_1)$, $\mf(z_2),\dots, \mf(z_n)$,
and $\mf(\infty)$ of $W$. 

Let $U$ be the open subset of $\Fl$ formed by
points $(V,E_{\bull})$ such that $E_{\bull}$ lies  in the
intersection of  $n+1$ open 
Schubert cells corresponding respectively to the $n+1$ distinguished
complete filtrations on $V$. This condition on $E_{\bull}$
is equivalent to the statement that for
$i=1,\dots,N+1$ and $z\in \{z_1,\dots,z_n\}$, 
the subspace $E_i$ has ramification sequence
$(a_{N+1-i+1}(z),\dots,$ $ a_{N+1}(z))$ at $z$.

The subset $U$ is dense in 
each fiber of $\Fl\to \Omega$. 

Consider the subset
  $\Fl^o\subseteq U$ formed by points $(V,E_{\bull})$ 
such that for $i=1,\dots,N-1$, the subspaces
$E_i\subset \cee[x]$ and $E_{i+1}\subset\cee[x]$ 
do not have common ramification points in $\cee-\{z_1,\dots,z_n\}$,  i.e.
their Wronskians do not have common roots in $\cee-\{z_1, \dots, z_n\}$.

\begin{lemma} $\Fl^o$ is open and dense in each fiber of $U\to \Omega$.
\end{lemma}
\begin{proof}
We recall the proof  from ~\cite{mukvar}, Lemma 5.19.  The
requirement for $E_i$ and $E_{i+1}$ to have common ramification at a given
$t\in\pone-\{z_1,\dots,z_s,\infty\}$ is at least ``$2$
conditions''. Taking into account the one parameter variation of $t$,
it is easy to see that $U-\Fl^o$ is of codimension at least one. Hence
the inclusions $\Fl^o\subseteq U\subseteq \Fl$ are open and dense.
\end{proof}

Suppose $(V, E_{\bull})\in \Fl^o$. 
Then for $i = 1, \dots , N+1$,  the polynomial
$\Wr(\, E_i )$ is divisible by $K_i$ and 
is of degree 
$i\,(d-i+1)\, - \,\sum_{\ell=1}^i  \, a_{N+1-i+\ell}(\infty). $ 
Introduce the polynomial $y_i$ by the condition
 $\Wr(\, E_i )\,  = \, K_i \, y_i$. The nonzero polynomial $y_i$
is defined up to multiplication by a number.

\begin{enumerate}
\item[($\alpha$)] If $l_i$ is the degree of $y_i$,  then
\begin{equation}\label{eqn2}
l_i\ = \ i\,(d-i+1)- \sum_{\ell=1}^i\ 
a_{N+1-i+\ell}(\infty)-\sum_{s=1}^n\sum_{\ell=1}^i a_{N+1-i+\ell}(z_s)\ .
\end{equation}
In particular, $y_{N+1}$ is of degree 0.
\item[($\beta$)] The polynomial
$y_i$ has no roots in the set $\{z_1,\dots, z_n\}$.
\end{enumerate}
The fact that $E_i$ and $E_{i+1}$ have no common ramification 
points in $\cee-\{z_1,\dots, z_n\}$ translates to the property
\begin{enumerate}
\item[($\gamma$)] The polynomials $y_i$ and $y_{i+1}$ have  no common roots.
\end{enumerate}
 Set $y_0 = 1$.

Suppose that $u_1, \dots, u_{N+1}$ is a basis of $V$ such that for
any $i=1, \dots, N$, the elements $u_1, \dots, u_i$ form a basis of $E_i$.
Let $Q_i \, =\,  \Wr(\,u_1, \dots , u_{i-1} , u_{i+1})\,/\, K_i$.  
\begin{lemma}\label{calcul}
$\Wr(\,y_i , Q_i)\, =\, T_i \, y_{i-1}\, y_{i+1}$ for $i=1,\dots,N$.
\end{lemma}
\begin{proof} 
$$
K_i^2\, \Wr(\,y_i,Q_i)\ =\ \Wr(\,K_i\,y_i,\,K_i\,Q_i)\ =\
\Wr(\,\Wr(\,u_1,\dots,u_i),\Wr(\,u_1 , \dots , u_{i-1} , u_{i+1}))\ .
$$
By  Lemma A.4 in ~\cite{mukvar}, the above quantity equals
$$\Wr(u_1,\dots,u_{i-1})\Wr(\,u_1 , \dots , u_{i+1} )\ =\
K_{i-1}\, y_{i-1}\, K_{i+1}\, y_{i+1}\ =\
T_i\, y_i\, y_{i+1}\, {K_i}^2\ .
$$
Finally, divide both sides by $K_i^2$.
\end{proof}

By Lemma ~\ref{calcul},  any multiple root of $y_i$ is
a root of either $T_i$, $y_{i-1}$, or $y_{i+1}$.
Clearly we have
\begin{enumerate}
\item[$(\delta)$]The polynomial $y_i$ has no multiple roots.
\item[$(\eta)$] There exist $\tilde{y}_i\in \cee[x]$ 
such that $\Wr(\,y_i , \tilde{y}_i)\, =\, T_i\,  y_{i-1}\, y_{i+1}$,\
namely, $\tilde{y}_i\, =\, Q_i$. 
\end{enumerate}
We translate condition $(\eta)$ into equations by using the following 

\begin{lemma}\label{new1}
Let $y\in\cee[x]$ be a polynomial with no multiple roots and 
$T\in\cee[x]$ any polynomial. Then  equation 
$\Wr(\, y , \tilde{y})\, =\, T$ has a solution $\tilde{y}\in\cee[x]$ 
if and only if $y$ divides $\Wr(\,y' , T )$.
\end{lemma}
This lemma follows from Lemma ~\ref{BA} below with $A=\cee$.

\bigskip

Consider         the space
$$
R\ =\ \prod_{i=1}^N \ \Bbb{P}(\cee_{ {l_i}}[x])
$$ 
where $l_i$ is given by (\ref{eqn2}).

Let $R^o$ be the open subset of $R$ formed by the tuples $(y_1,\dots,y_N)$ 
 satisfying conditions $(\alpha)-(\delta)$. Let $\mathcal{A}$ be the 
subset of $R^o$  defined by the condition 
\begin{equation}\label{baa}
y_i 
\quad
{\rm divides}
\ \
 \Wr(\,y_i', \, T_i y_{i-1}y_{i+1})
\ \
{\rm for}\ \ i = 1, \dots ,  N\ .
\end{equation}
Using the monicity of $y_i$ and long division, we can write the divisibility 
condition as a system of equations in the coefficients of $y_i$, $y_{i-1}$,
 and $y_{i+1}$. Hence $\mathcal{A}$ is a closed subscheme of $R^o$.

Consider the morphism 
$$
\Theta\ :\ \Fl^o\ \to\ R^o\ , 
\qquad
(V,E_{\bull})\ \mapsto\ (y_1,\dots,y_N)\
=\ (\,\Wr(\, E_1 ) /  K_1, \dots, \Wr(\, E_N ) / K_N)\ .
$$ 
For $x\in \Fl^o$, condition $(\eta)$ holds and by Lemma ~\ref{BA}, 
$\Theta$ induces a morphism of schemes $\Theta:\Fl^o\to\ma$.

\begin{theorem}\label{mainte}
The morphism \ $\Theta\,:\,\Fl^o\,\to\, \ma$\
is an isomorphism of schemes.
\end{theorem}
It is proved in ~\cite{mukvar} 
that $\Theta$ is a bijection of sets. 
In Section ~\ref{proofe} we will extend the argument 
of ~\cite{mukvar} to prove Theorem ~\ref{mainte}. 

\bigskip

The following corollaries of Theorem ~\ref{mainte} 
use the notion of the geometric multiplicity of an 
irreducible scheme. 
This notion, as well as its relation to intersection 
theory is reviewed in the appendix.    

\begin{corollary}\label{MLK}
Let $x\in \Omega$. Let $m(x)$ be the length of the local ring of 
$\Omega$ at $x$. Let $C$ be the irreducible component of $\ma$ 
which, as a point set, is 
$\Theta(\pi^{-1}(x)\cap \Fl^o)$,  see the Cartesian square 
(\ref{cartsquare}). 
Then, the  geometric multiplicity of $C$ equals $m(x)$.
In particular,  the geometric multiplicity of $C$ equals  the geometric
multiplicity of $\Omega$ at $x$.
\end{corollary}
\begin{proof}
Let $\Omega_x$ be the component of 
$\Omega $ containing $x$. As a set,   
$\Omega_x$ is just the point $x$. 
Consider the irreducible component $\mi$ of $\Fl^o$ 
containing $\pi^{-1}(x)$, 
\begin{equation}
\xymatrix{\mi\ar[r]\ar[d] & \Fl\ar[r]\ar[d]^{\pi}  &                \Fl(\mv)\ar[d]\\
           \Omega_x \ar[r] & \Omega\ar[r]                          & \Gr(N+1,W) }
\end{equation}
Since $\pi$ is a locally trivial fiber bundle, the morphism $\mi\to
\Omega_x$ is a fiber bundle with smooth fiber for the Zariski topology
on the scheme $\Omega_x$ . The multiplicity of $\mi$ is the same as
the multiplicity of its dense subset $\mi\cap \Fl^o$. The corollary
now follows from the theorem and Proposition ~\ref{pga}.
\end{proof}

\begin{corollary} We have an equality of cohomology classes
in $ H^*(\Gr(N+1,W))$, 
$$
[\Omega_{\aA(\infty)} (\mf(\infty))]\ \cdot\
\prod_{s=1}^n \
[\Omega_{\aA(z_s)}(\mf(z_s))]\ =\ c\
 [\text{ class of a point }]
$$ 
where $c$ is the sum of the geometric multiplicities of  
irreducible components of $\mathcal{A}$.
\end{corollary}
\subsection{Critical point equations}
Consider the space $\tilde R = \prod_{i=1}^N \cee^{l_i}$ with
coordinates $(t^{(i)}_j)$, where $i=1,\dots,N$, \ $j=1,\dots,l_i$.
The product of symmetric groups 
\linebreak
$\Sigma = \Sigma_{l_1}\times \dots\times \Sigma_{l_N}$
acts on $\tilde R$ by permuting coordinates with the same upper index.
Define a map
$$
\Gamma\ :\ \tilde R \ \to\ R\ , 
\qquad
(t^{(i)}_j)\ \mapsto \ (y_1,\dots,y_N)\ ,
$$  
where
$
y_i\ =\ \prod_{j=1}^{l_i}\ (x-t^{(i)}_j)\ .
$
Define the scheme $\tilde \ma$ by the condition $\tilde \ma =
 \Gamma^{-1}(\ma)$. The natural map $\tilde{\ma}\to\ma$ is finite and
 \'{e}tale.  The scheme $\tilde{\ma}$ is $\Sigma$-invariant. The
 scheme $\tilde{\ma}$ lies in the $\Sigma$-invariant subspace $\tilde
 R^o$ of all $(t^{(i)}_j)$ with the following properties for every $i$:
\begin{enumerate}
\item[$\bullet$] The numbers $t^{(i)}_1, \dots , t^{(i)}_{l_i}$ are distinct;
\item[$\bullet$]
The sets $\{t^{(i)}_1, \dots , t^{(i)}_{l_i}\}$
and $\{t^{(i+1)}_1, \dots , t^{(i+1)}_{l_{i+1}}\}$ do not intersect;
\item[$\bullet$]
The sets $\{t^{(i)}_1, \dots , t^{(i)}_{l_i}\}$
and $\{z_1, \dots , z_n\}$ do not intersect.
\end{enumerate}
\begin{lemma}\label{oldone}${}$

\begin{enumerate}
\item[$\bullet$]
 The connected components of $\ma$ and $\tilde{\ma}$ are irreducible.
\item[$\bullet$]
  The reduced schemes underlying $\ma$ and $\tilde{\ma}$ are smooth. 
\item[$\bullet$]
 If $C$ is a connected component of $\ma$, then the group
$\Sigma$ acts transitively on the connected components of
$\Gamma^{-1}(C)$.
\end{enumerate}
\end{lemma}
\begin{proof}
By Theorem ~\ref{mainte}, $\ma$ is  isomorphic to $\Fl^o$. The
reduced scheme underlying $\Fl^o$ is smooth. Therefore the reduced
scheme underlying $\ma$ is smooth. Since $\Gamma$ is \'{e}tale, we
deduce that the reduced scheme underlying $\tilde{\ma}$ is also
smooth.

The smoothness conclusions immediately imply the irreducibility of
connected components of $\ma$ and $\tilde{\ma}$.

The transitivity assertion follows from the fact that $\Gamma$ is a
Galois covering with Galois group $\Sigma$.
\end{proof}

Consider on $\tilde R^o$ the regular rational function 
$$
\Phi(t^{(i)}_j)\ =\ \prod_{i=1}^N \prod_{j=1}^{l_i}  T_i(t^{(i)}_j)^{-1} \
\prod_{i=1}^{N-1} \prod_{j=1}^{l_i}\prod_{k=1}^{l_{i+1}}
(t^{(i)}_j- t^{(i+1)}_k)^{-1}\
\prod_{i=1}^{N}\prod_{1\leq j< k \leq l_i}\!
(t^{(i)}_j-t^{(i)}_k)^{2}\ .
$$ 
This $\Sigma$-invariant function is called {\it the master function
associated with the basic situation}.

Define the scheme $\tilde \ma'$ as the subscheme in $\tilde R^o$ of critical points
of the master function.

\begin{lemma}\label{bethe eqn lemma}
The subschemes $\tilde \ma$ and $\tilde \ma'$ of  $\tilde R^o$ coincide.
\end{lemma}
\begin{proof} 
The subscheme $\tilde{\ma}$ is
defined by divisibility conditions
(\ref{baa}).
By Lemma ~\ref{equationes}, the divisibility condition
(\ref{baa}) for a fixed $i$, reduces to the critical point equations
for the function 
$$ 
\prod_{j=1}^{l_i}  T_i(t^{(i)}_j)^{-1}
\prod_{j=1}^{l_i}\prod_{k=1}^{l_{i+1}} (t^{(i)}_j- t^{(i+1)}_k)^{-1}
\prod_{j=1}^{l_i}\prod_{k=1}^{l_{i-1}}(t^{(i-1)}_k- t^{(i)}_j)^{-1}\!
\prod_{1\leq j< k \leq l_i} \!(t^{(i)}_j-t^{(i)}_k)^{2}\ 
$$ 
of variables $t^{(i)}_1, \dots , t^{(i)}_{l_i}$. Then the system of
divisibility conditions (\ref{baa}) for all $i$ together is just the
critical scheme of $\Phi$. This concludes the proof.
\end{proof}

Let $(V,E_{\bull})\in \Fl^o$. Denote $\bs y = \Theta(V,E_{\bull})\in\ma$.
Pick a point $\bs t \in \Gamma^{-1}(\bs y)$. Let $C$ be the unique
irreducible component of $\ma$ containing $\bs y$ and $\tilde{C}$ the
unique irreducible component of $\tilde{\ma}$ containing
$\bs t$.
\begin{theorem}\label{crit point multiplicity}
The geometric multiplicity of the scheme $\Omega$ at $V$ equals
the geometric multiplicity of $\tilde{C}$.
\end{theorem}
\begin{proof}
The morphism $\tilde{C}\to C$ is \'{e}tale. By Proposition
~\ref{pga}, the geometric multiplicity of $\tilde{C}$ coincides with
that of $C$. Now the theorem follows from Corollary ~\ref{MLK}.
\end{proof}

The group $\Sigma$ acts on the set of connected components of
$\tilde{\ma}$.  For an orbit of this action,  define its
geometric multiplicity to be the geometric multiplicity
of any member of the orbit. 

\begin{corollary} We have an equality of cohomology classes
in $ H^*(\Gr(N+1,W))$, 
$$
[\Omega_{\aA(\infty)} (\mf(\infty))]\ \cdot\
\prod_{s=1}^n \
[\Omega_{\aA(z_s)}(\mf(z_s))]\ =\ c\
 [\text{ class of a point }]
$$ 
where $c$ is the number of orbits for the action of $\Sigma$ on the
connected components of $\tilde{\ma}$ counted with geometric
multiplicity.
\end{corollary}

\section{Proof of Theorem ~\ref{mainte}}\label{proofe}

\subsection{Admissible modules}

Let $(A,\emm)$ be a local ring with residue field $\cee$. A submodule $V\subset A[x]$ is said 
to be {\it an admissible submodule of rank $N+1$}, if
\begin{enumerate}
\item [$\bullet$]
The submodule $V\subset A[x]$  is a free $A$-module of rank $N+1$;
\item [$\bullet$]
The morphism $V\tensor {A/\emm}\to ({A/\emm})[x]=\cee[x]$ is injective.
\end{enumerate}
An admissible submodule $V\subset A[x]$  
is said to have {\it the ramification sequence $\aA(z)=(a_1,\dots,a_{N+1})$ 
at $z\in\cee$} \ {} if\ {} $V$\  has an $A$-basis 
$$
(x-z)^{N+1-1+a_1}f_1,\ (x-z)^{N+1-2+a_2}f_2 ,\ \dots ,\
 (x-z)^{N+1-(N+1)+a_{N+1}}f_{N+1}
$$
with $f_i(z)\in  A^*$ for $i = 1, \dots , N+1$. 

An admissible submodule $V\subset A[x]$  
is said to have {\it the ramification sequence $\aA(\infty) = (a_1,\dots,a_{N+1})$ 
at $\infty$}\ {} if\ {}
$V$ has an $A$-basis $f_1,$ $f_2 ,$ $ \dots $, $f_{N+1} $ with
monic $f_i$  and $\deg\,f_i\, = \linebreak
d-(N+1)+i-a_i$ for $i=1,\dots,N+1$. 

\begin{remark}
Admissible modules $V\subset A[x]$ of rank $N+1$ are in one-to-one
correspondence with morphisms $\spec(A)\to \Gr(N+1,\cee[x])$. Intuitively,
if $A=\cee[[t]]$, this is a formal holomorphic map of a $1-$disc into
$\Gr(N+1,\cee[x])$.

An admissible submodule $V\subset A[x]$ may not have a ramification
sequence at a given $z\in\cee\cup\{\infty\}$. This corresponds
intuitively, to the case when the formal map considered above does not
remain in a Schubert cell. If $V$ has a ramification sequence at $z$,
then the ramification sequence is unique (equal to the ramification
sequence of the subspace $V\tensor (A/\emm)\subset \cee[x]$ at $z$).
\end{remark}

For $f \in A[x]$, we denote by $\bar{f}$ the corresponding polynomial in
$({A/\emm})[x] = \C[x]$. 

The following standard lemma is proved in Section ~\ref{setupproof}.

\begin{lemma}\label{setup}
Let $A$ be a local ring with residue field $\C$ and $V\subset A[x]$ a submodule. 
Then the following statements are equivalent.
\begin{enumerate}
\item The submodule $V$ is admissible.
\item The submodule $V$  
is a finitely generated $A$-module, and there is an
$A$-module decomposition $A[x] \,=\,V\,\oplus\, M$ for some $A$-module $M$.
\item For some $k$, there exist $u_1, \dots , u_k \in V$ such that $V$ is the $A$-span
of $u_1, \dots , u_k$ and the elements $\bar{u}_1,\dots,\bar{u}_k \in
\cee[x]$ are linearly independent over $\C$.
\end{enumerate} 
\end{lemma}

\begin{lemma}\label{closed}
Let $A$ be a local ring with residue field $\cee$  and $V\subset A[x]$
 an admissible submodule of rank  $N+1$. 
Let $u_1,\dots,u_{N+1}$ be a basis of $V$ as an $A$-module. Suppose
$v\in A[x]$ satisfies the equation 
$\Wr(\,u_1, \dots , u_{N+1}, v)\,=\, 0$. Then $v\in V$.
\end{lemma}

Lemma \ref{closed} is proved in Section ~\ref{remainder}.

\subsection{Proof of Theorem ~\ref{mainte}}

Our first objective will be to show that $\Theta\,:\,\Fl^o\to R^o$ is a closed
embedding of schemes. 
Our second objective will be to show that $\Fl^o\to \ma$ is an isomorphism
of schemes. 

Clearly  the morphism $\Theta\,:\,\Fl^o\to R^o$
extends to a morphism of projective schemes $\tiTheta \,:\, \Fl \to R$ 
and $\Fl^o\,=\,\tiTheta^{-1}(R^o)$. The morphism $\tiTheta$ is closed as
 a morphism of projective schemes.

We  achieve the first objective by showing that $\tiTheta$ is an
embedding. By Lemma ~\ref{schema}, it will be enough to show  that for
every local ring $A$ over $\cee$, the morphism $\Fl(A)\to R(A)$ is an
injective map of sets. 

To achieve the second objective, by Lemma ~\ref{schema2}, 
it will be  enough to show that for any local ring $A$ over $\C$, the
induced map $\Fl^o(A)\to\ma(A)$ is a set theoretic surjection.

To use these criteria, we need to define the sets $\Fl(A)$,\ $\Fl^o(A)$,\ $R(A)$, and
$\ma(A)$ more explicitly.

\medskip

By Proposition ~\ref{functor}, the set $\Fl(A)$ is the set of
pairs $(V,E_{\bull})$, \
such that 
\begin{enumerate}
\item[$\bullet$]
$V\subset A[x]$ is an admissible submodule
of rank $N+1$;
\item[$\bullet$]
$E_{\bull}=(E_1\subsetneq
E_2\subsetneq\dots\subsetneq E_{N+1} = V)$ is a filtration by
admissible submodules;
\item[$\bullet$]
$V$ has the ramification sequence 
$(a_{1}(z),\dots, a_{N+1}(z))$ at each $z\in
\{z_1,\dots,z_n,\infty\}$.
\end{enumerate}
The subset $\Fl^o(A)$ is the set of
$(V,E_{\bull})\in \Fl(A)$ such that the induced point in
$\Fl(\cee)$ is a point of $\Fl^o$.

The set $R(A)$ is the set 
$$
\{ (y_1,\dots,y_N) \in \ \prod_{i=1}^N A_{l_i}[x]\ \ | \ \  
\bar{y}_i \neq 0 ,\  i = 1 , \dots , N \}
$$
 modulo the equivalence relation
 $(y_1,\dots,y_N) \sim (\tilde{y}_1,\dots,\tilde{y}_N)$ if there exist
 $a_i\in A^*$ such that $a_i\, y_i\, = \,\tilde{y}_i$\ for\ $i = 1, \dots , N$.

The subset $\ma(A)$ consists of elements
$(y_1,\dots,y_N)\in R(A)$ such that
\begin{enumerate}
\item[$\bullet$]
 For $i=1,\dots,N$, condition $(\eta)$ holds. By Lemma ~\ref{BA},
this implies that the polynomial 
$y_i$ divides $\Wr(\,y_i',\, T_i\, y_{i-1}\,y_{i+1})$. Here  $y_0\,=\,
y_{N+1}\,=\,1$.
\item[$\bullet$]
 The reduction $(\bar{y}_1,\dots,\bar{y}_{N})\in R(\cee)$ is a point of $R^o$.
\end{enumerate}

\begin{lemma}\label{settheory} The morphism
$\tiTheta:\Fl\to R$ is a closed embedding of schemes.
\end{lemma}
\begin{proof}
We need to show that $\tiTheta:\Fl(A)\to R(A)$ is a set theoretic
injection for any local ring $A$.

Suppose $(V,E_{\bull})$ and $(V',E'_{\bull})$ are two points of
$\Fl(A)$ with $\tiTheta(V,E_{\bull})=\tiTheta(V',E'_{\bull})$. Pick
bases $(u_1,\dots,u_{N+1})$ and $(v_1,\dots,v_{N+1})$ for $V$ and $V'$
respectively so that for all $i$, the admissible submodule
$E_i$ is the $A$-span of $u_1,\dots, u_i$ and 
the admissible submodule $E'_i$ is the $A$-span of $v_1,\dots, v_i$.

The hypothesis implies that $\Wr(\,u_1,\dots,u_i)\,=\,c_i\,
\Wr(\,v_1,\dots,v_i)$
with $c_i \in A^*$. Clearly
$E_1=E'_1$. Assume by induction that $E_i = E_i'$. Then
$$
\Wr(\,v_1,\dots,v_{i+1})\ =\ c\ 
\Wr(\,u_1,\dots,u_{i},v_{i+1})\ =\ c'\ \Wr(\,u_1,\dots,u_{i},u_{i+1})
$$
for $c, c' \in A^*$. Therefore 
$$
\Wr(\,u_1,\dots,u_{i},\,c v_{i+1} - c' u_{i+1})\ =\ 0\ .
$$
Lemma ~\ref{closed} implies $c v_{i+1} -c' u_{i+1} \in E_{i}$ 
and hence $E_{i+1} = E'_{i+1}$. 
\end{proof}

We need to show that for any local ring $A$ over $\C$, the
induced map $\Fl^o(A)\to\ma(A)$ is a set theoretic surjection.
But this claim on surjectivity follows from the following
theorem on the existence of solutions to Wronskian equations.

\medskip

Let $T_0, \, T_1,\,\dots,T_N \,\in\cee[x]\subseteq A[x]$ be non-zero
polynomials. Let $S \,\subset \,\cee$ be the union of their zero sets.  Set
$K_i\,=\,T_0^i\, T_1^{i-1}\,T_2^{i-2}\,\cdots \,T_{i-1}$ for $i=1,\dots, N+1$.

Let $y_1,\dots, y_N\in A[x]$ be monic polynomials of
arbitrary degree.  Set $y_0 \,=\, y_{N+1}\, =\, 1$. 

For $z\in S$ and $i = 1,\dots,N+1$, define 
$$
e_{i}(z) \, =\, i-1 +  \sum_{j=0}^{i-1}\, \ord_z  T_j\ ,
\qquad
c_{i}\, =\, i-1 + \deg\,y_i\, -\, 
\deg\,y_{i-1}\, +\, \sum_{j=0}^{i-1}\, \deg\,T_j\ .
$$

\begin{theorem}\label{Main} Under these conditions assume that
for all $i$
\begin{enumerate}
\item[$\bullet$] The polynomial $ \bar{y}_i\in \C[x]$ has no
multiple roots,  no roots in $S$,  and is coprime to  $\bar{y}_{i+1}$;
\item[$\bullet$] The polynomial
$y_i$ divides $\Wr(\,y_i',\,T_i\, y_{i-1}\,y_{i+1})$.
\end{enumerate}
Then, there exist $u_1,\dots,u_{N+1}\in A[x]$ with the following properties.
Set $E_i$ to be the $A$-span of $u_1,$ $\dots ,$ $u_i$ for $i = 1, \dots , N+1$ 
and set  $V\,=\,E_{N+1}$. Then for all $i$
\begin{enumerate}
\item  $\Wr(\,u_1, \dots , u_i )\ =\, K_i\, y_i$;
\item $E_i \subset A[x]$ is admissible;
\item $E_i$ has ramification sequences at each $z\in
S\cup\{\infty\}$.  The set of exponents at $z\in S$ is
$\{e_1(z),\dots,e_{i}(z)\}$. The set of exponents at $\infty$ is
$\{c_1,\dots,c_{i}\}$.
\end{enumerate}
\end{theorem}

Theorem ~\ref{Main} is proved in Section ~\ref{proof of Main}. 

Let $(y_1,\dots,y_N)\in \ma(A)$. Apply Theorem \ref{Main} and
obtain a point $(V,E_{\bull})$. It is easy to see that
$(V,E_{\bull})\in \Fl^o(A)$ and
$\Theta(V,E_{\bull})=(y_1,\dots,y_N)$. The proof of Theorem
~\ref{mainte} is complete.
\subsection{Proof of Lemma ~\ref{closed}}\label{remainder}
Let $\emm \subset A$ be the maximal ideal.
Let $u_1, \dots , u_{N+1}\,\in\, V$ be a basis, \ 
$u_i\, = \,\sum_{l}\,a^i_{l}\, x^l$ with $a^l_i \in A$.

If necessary changing the basis, 
we may assume the basis has the following property.
There exist nonnegative integers
$k_1 > \dots > k_{N+1}$ such that
$$
a^j_{k_i}\,=\,0 \ {\rm for}\, j\neq i\ ; \ {}
\qquad
a^i_{k_i}\,\in\, A-m\ ; \ {}
\qquad
a^i_{j}\,\in m \ {\rm for}\ j\,>\,a^i_{k_i}\ .
$$
Let $v\,=\,\sum_{l}\,b_l\,x^l$ be a nonzero polynomial such that
$\Wr(\,u_1, \dots , u_{N+1}, v)\,=\, 0$. We may assume that $b_{k_i} =
0$ for all $i$. We shall prove that this leads to contradiction.

Recall that for any nonzero $a\in A$, there is a unique $r$ such that
 $a\in m^r-m^{r+1}$, see Krull's intersection theorem (
 ~\cite{mat}, Theorem 8.10).

Let $r$ be the smallest number such that some $b_l$ is in
$m^r-m^{r+1}$, and $p$ the largest \linebreak
index such that $b_p \in
m^r-m^{r+1}$. Clearly $p  \not\in \{k_1 , \dots , k_{N+1}\}$. The polynomial
\linebreak
$\Wr(\,u_1 , \dots , u_{N+1} , v)$ is nonzero since in its decomposition
into monomials, the coefficient of the monomial $\Wr(\,x^{k_1} , \dots
, x^{k_{N+1}} , x^p )$ belongs to $m^r - m^{r+1}$.

\subsection{Proof of Theorem ~\ref{Main}}\label{proof of Main}
The proof follows ~\cite{mukvar}. 
Call a tuple $(y_1,\dots,y_N)$ {\it fertile} if it
satisfies the conditions of Theorem ~\ref{Main}.

For  $i = 1 , \dots , N$, define the process of
reproduction in the $i$-th direction. Namely, find a solution 
$\tilde{y}_{i} \in A[x]$ to the equation 
$\Wr(\,y_{i} , \tilde{y}_{i} )\,=\, T_{i}\,y_{i+1}\,y_{i-1}$.  
If \ $\tilde{y}_{i}$ is a solution, then for  $c\in A$, the polynomial
$\tilde{y}_{i} + c y_i$ is a solution too. Add to $\tilde y_i$ the term
$c y_i$ if necessary, and obtain  a monic ${\tilde {y}}_i$ such that
its reduction  modulo the
maximal ideal does not have roots in $S$,
does not have multiple roots, and has no common roots with reductions of
$y_{i-1}$ or $y_{i+1}$. The transformation from the tuple $(y_1, \dots , y_n)$ to  
$(y_1, \dots , y_{i-1}, \tilde y_i , y_{i+1}, \dots , y_N)$
is the process of reproduction in the $i$-th direction. 

\medskip

{\bf Claim:} {\it The tuple $(y_1,\dots,\tilde{y}_i,\dots, y_N)$ is fertile.} 

\medskip

To prove the claim it is enough to show that
\begin{equation}\label{questionmark1}
y_{i-1} \ \
 {\rm divides}
\ \
\Wr(\,y_{i-1}' , T_{i-1}\,y_{i-2}\,\tilde{y}_{i})
\end{equation}
and 
\begin{equation}\label{questionmark2}
y_{i+1}\ \
{\rm divides}
\ \ 
 \Wr(\,y_{i+1}' , \,T_{i+1}\,y_{i+2}\,\tilde{y}_{i})\ .
\end{equation}
We prove (\ref{questionmark1}). Statement (\ref{questionmark2}) is
proved similarly. 

\medskip

Clearly, $y_{i-1}$ divides $\Wr(\,y_{i} , \ \tilde{y}_{i})$. By
assumption, $y_{i-1}$ divides $\Wr(\,y_{i-1}' ,
\,T_{i-1}\,y_{i-2}\,y_{i})$. By Jacobi's rule,
$$
\Wr(\,y_{i-1}' , \,T_{i-1}\ y_{i-2}\,
\tilde{y}_{i})\ y_{i}\ -\ \Wr(\,y_{i-1}' ,\,T_{i-1}\,y_{i-2}\,
{y}_{i} )\ \tilde{y}_{i}
\ =\ \Wr(\,y_{i} ,,\tilde{y}_{i})\ T_{i-1}\,
y_{i-2}\,y_{i-1}'\ .
$$ 
Hence  $y_{i-1}$ divides
$\Wr(\,y_{i-1}' , \,T_{i-1}\,y_{i-2}\,\tilde{y}_{{i}})\ y_{i}$.
By assumption, the ideal $(\bar{y}_{i-1}, \bar{y}_i)\ =\ \cee[x]$.
By  Lemma ~\ref{mon1}, this implies that the ideal $(y_{i-1},y_i)=A[x]$. 
Now use Lemma ~\ref{mon2} to see that $y_{i-1}$ divides
$\Wr(\,y_{i-1}' , \,T_{i-1}\,y_{i-2}\,\tilde{y}_{{i}})$.
This proves (\ref{questionmark1}) and the claim.

\medskip

To construct the polynomials $u_1$, $\dots,$ $u_{N+1}$ we do the
following. We set $u_1 \,=\, K_1\,  y_1$. To 
construct the polynomial $u_{i+1}$, for
$i = 1, \dots , N$, 
we perform the simple reproduction procedure
in the $i$-th direction and obtain the tuple
$(y_1, \dots , \tilde y_i, \dots , y_N)$. Next we perform for   
$(y_1, \dots , \tilde y_i, \dots , y_N)$ the simple reproduction procedure in the
direction of $i-1$ and obtain 
$(y_1, \dots , \tilde y_{i-1} , \tilde y_i , \dots , y_N)$. 
We repeat this procedure all the way until the simple reproduction procedure in the
1st direction and obtain 
$(\tilde y_1, \dots , \tilde y_{i-1} , \tilde y_i , \dots , y_N)$. 
We set $u_{i+1} = K_1 \tilde y_1$.

\medskip

{\bf Claim:} {\it We have
$\Wr(\,u_1 , \dots , u_i)\ =\ K_i\,y_i$ \  for $i = 2, \dots , N+1$.}

\medskip
We prove the claim by 
induction. Let $(y_1,\dots,\tilde{y}_{i},\dots,y_N)$ be the tuple
obtained by the simple reproduction in the $i$-th direction. Apply induction to the
tuple $(y_1,\dots,\tilde{y}_{i},\dots,y_N)$ to obtain
$\Wr(\,u_1 , \dots , u_{i-1} , u_{i+1})\,=\,K_i\,\tilde{y}_i.$ 
Induction applied to the tuple $(y_1,\dots,{y}_{i},\dots,y_N)$ gives 
$\Wr(\,u_1 , \dots , u_{i-1} , u_{i} )\ =K_i \ y_i$. 
Now
$$
\Wr(\,u_1 , \dots , u_{i+1})
\Wr(\,u_1 , \dots , u_{i-1} )\ =\ \Wr(\,\Wr(\,
u_1 , \dots , u_{i-1} , u_i ),\,\Wr(\,u_1 , 
\dots , u_{i-1} , u_{i+1}))
$$
$$
=\,
\Wr(\,K_i\,y_i , \,K_i\, \tilde{y_i})\ =\ K_i^2\ T_i\ y_{i_1}\ y_{i+1}\ .
$$ 
By
induction we also have $\Wr(\,u_1 , \dots , u_{i-1})$ $ =$
$ K_{i-1}\,y_{i-1}$
and $T_i\,=\,K_{i+1}\,K_{i-1} / K_i^2.$
The above equation rearranges to 
$\Wr(\,u_1 , \dots , u_{i+1})\, K_{i-1}\, y_{i-1} = K_{i+1}\,
K_{i-1}\, y_{i-1}\, y_{i+1}$, 
which gives the desired equality
$\Wr(\,u_1 , \dots , u_{i+1} ) = K_{i+1}\, y_{i+1}$.
\medskip

Return to the proof of Theorem ~\ref{Main}. 
Let $E_i$ be the $A$-span of $u_1,\dots,u_{i}$. Then
$\Wr(\,\bar{u}_1 , \dots , \bar{u}_{i} )
 = \bar{K}_{i}\, \bar{y}_{i} = K_{i}\, \bar{y}_{i}\, \neq 0$.
By part (3) of Lemma ~\ref{setup}, the submodule
$E_i \subset A[x]$ is admissible.
This proves (1) and (2) of the theorem.

Now we will calculate the exponents of $E_i$ at $z\in S$. The
exponents of $E_i$ at $\infty$ are calculated similarly.

By induction assume that the set of exponents of $E_i$ at $z$ is
 $\{e_1(z),\dots,e_i(z)\}$. That means that $E_i$ has an $A$-basis $v_1 , \dots , v_i$
such that $v_j\,=\,(x-z)^{e_j(z)}\,a_j$ with $a_j(z)\in A^*$ for all $j$. Hence
$\Wr(\,E_i)$  is of the form $(x-z)^{b}\, a$  with
$$
b\ =\ \sum_{\ell=1}^i\,e_{\ell}(z)\ -\ \frac{i\,(i-1)}2 \ ,
\qquad a(z)\in A^*\ .
$$ 
We have  $b = \ord_z\,K_i$ since $\Wr(\,E_i) = K_i \,y_i$.

Let $v_{i+1}\in E_{i+1}$ be such that $v_1, \dots , v_i, v_{i+1}$ is
an $A$-basis of $E_{i+1}$. We may assume that $ v_{i+1} =
\sum_{\ell=0}^{d}\, c_{\ell}\ (x-z)^\ell\ $ with coefficients
$c_{\ell}$ equal to zero for $\ell\,\in\,
\{e_1(z),\dots,e_i(z)\}$. Let $p$
be the smallest integer with $c_p\,\neq \,0$. It is easy to see that the
Wronskian of $E_{i+1}$ equals $(x-z)^{b\,+\,p\,-\,i}\, h$ where $h(z)$ is
$c_p$ up to multiplication by an element in $A^*$.  From equation
$\Wr(\,E_{i+1})\,=\,K_{i+1}\,y_{i+1}$ we deduce that $\Wr(\,E_{i+1})$ is
of the form $(x-z)^{\ord_z\,K_{i+1}}\, g$ with $g(z)\,\in\,
A^*$.  Hence $b+p-i=\ord_z\,K_{i+1}$, \ $p=e_{i+1}(z)$, and $c_p\,\in\,
A^*$. 

This reason shows that
$e_{i+1}(z)$ does not belong to the set
$\{e_1(z),\dots, e_i(z)\}$. Thus $E_{i+1}$
has a ramification sequence at $z$ and the set of exponents at $z$ is
$\{e_1(z),\dots, e_i(z), e_{i+1}(z)\}$. The induction step is
complete.

\section{Schubert  cells and critical points}\label{extendedbasic}

Let $\Sigma_{N+1}$ be the permutation group of the set $\{1, \dots , N+1\}$
and $w\in \Sigma_{N+1}$. Assume that a basic situation
of Section  ~\ref{secion} is given. 
In Section ~\ref{secion}, we defined a flag bundle $\Fl$ over
$\Omega$. We also observed that $V\in\Omega$ has $n+1$
distinguished complete flags on $V$ induced from the 
complete flags $\mf(z_1),\dots,\mf(z_n)$ and $\mf(\infty)$ of $W$.

Let $\{c_1>\dots>c_{N+1}\}$ be the set
$\{d-(N+1)+j-a_j(\infty)\ |\ j=1,\dots,N+1\}$ of exponents of $V$ at
$\infty$. 

We define a subset $\Fl^o_w \subset \Fl$ as follows. Let $U_w$ be the 
subset of $\Fl$ formed by points $(V,E_{\bull})$ such that 

\begin{itemize}
\item $E_{\bull}\in \Fl(V)$ lies in the intersection of $n$ open
Schubert cells corresponding respectively
 to the $n$ distinguished complete
flags on $V$ induced from the flags $\mf(z_1), \dots , \mf(z_n)$.
This condition on $E_{\bull}$
is equivalent to the statement that for
$i=1,\dots,N+1$ and $z\in \{z_1,\dots,z_n\}$, 
the subspace $E_i$ has ramification sequence
$(a_{N+1-i+1}(z),\dots,$ $ a_{N+1}(z))$.

\item $E_{\bull}\in \Fl(V)$ lies in the Schubert cell, corresponding to
the permutation $w$ and the distinguished complete flag on $V$ induced
from the flag 
$\mf(\infty)$. This condition on $E_{\bull}$
is equivalent to the statement that the set of 
exponents of $E_i$ at $\infty$ is 
$\{c_{w(1)},\dots,c_{w(i)}\}$.

\end{itemize}
We define   $\Fl^o_w\subseteq U_w$ as the subset of points $(V,E_{\bull})$ 
such that for all $i$, the subspaces
$E_i\subset \cee[x]$ and $E_{i+1}\subset\cee[x]$ 
do not have common ramification points in $\cee-\{z_1,\dots,z_n\}$.

Notice that the subset  $\Fl^o_w$ may be empty if $w$ is not the identity element in
$\Sigma_{N+1}$.

\begin{lemma}
The morphism $\Fl^o_w\to \Omega$ is  smooth.
\end{lemma}
\begin{proof}
Let $\mathcal{J}$ be the 
subset of $\Fl$ formed by points
$(V,E_{\bull})$ such that $E_{\bull}$ lies in the Schubert cell
corresponding to the permutation $w$ and the distinguished flag
on $V$ induced from $\mf(\infty)$. It is easy to see that $U_w$ is an open subset of
$\mathcal{J}$ and therefore it suffices to show that $\mathcal{J}\to
\Omega$ is smooth. \

Denote by $\mv_{\aA(\infty)}$ the pull back of $\mv$ to
$\Omega^o_{\aA(\infty)}(\mf(\infty))\hookrightarrow\Gr(N+1,W)$.  There
is a distinguished section of $ \Fl(\mv_{\aA(\infty)})\ \to\
\Omega^o_{\aA(\infty)}(\mf(\infty)) $, \, see
Section \ref{1.1}.
 Let $G_w^{\mf(\infty)} \subset
\Fl(\mv_{\aA(\infty)})$ be the
part corresponding to $w$ in the partition of $\Fl(\mv_{\aA(\infty)})$
into Schubert cells associated to this distinguished section.

There is a fiber square 
$$
\xymatrix{\mathcal{J}\ar[r]\ar[d]   &  G_w^{\mf(\infty)}\ar[d]^{p}
\\
 \Omega \ar[r]  & \Omega^o_{\aA(\infty)}(\mf(\infty)) }
$$
The morphism $p$ is clearly smooth. The morphism
$\mathcal{J}\to \Omega$ is the base change of a smooth morphism and hence
it is smooth too.
\end{proof}
Let $(V,E_{\bull})\in \Fl^o_w$. 
For $i = 1, \dots , N+1$,  the polynomial
$\Wr( \,E_i )$ is divisible by $K_i$ and 
has degree 
$ \sum_{j=1}^i\,c_{w(j)}\,-\,i(i-1)/2$. 
Introduce the polynomial $y_i$ by the condition
 $\Wr(\, E_i )\,  = \, K_i \, y_i$. 

\begin{enumerate}
\item[$(\alpha)_w$] If $l^w_i$ is the degree of $y_i$, then
$l^w_i\,=\,\sum_{j=1}^i \,c_{w(j)}\,-\,i(i+1)/2\,-\,\deg\,K_i$. In particular,
$y_{N+1}$ is of degree 0.
\item[$(\beta)_w$] The polynomial
$y_i$ has no roots in the set $\{z_1,\dots, z_n\}$.
\item[$(\gamma)_w$] The polynomials $y_i$ and $y_{i+1}$ have  no common roots.
\item[$(\delta)_w$]The polynomial $y_i$ has no multiple roots.
\item[$(\eta)_w$] There exist $\tilde{y}_i\in \cee[x]$ 
such that $\Wr(\,y_i , \tilde{y}_i)\, =\, T_i\,  y_{i-1}\, y_{i+1}$.
\end{enumerate}

Consider  the space
$$
R_w\ =\ \prod_{i=1}^N \ \Bbb{P}(\cee_{ {l^w_i}}[x])
$$ 
where $l^w_i$ is given by property $(\alpha)_w$. 

Let $R_w^o$ be the open subset of $R_w$ formed by the tuples $(y_1,\dots,y_N)$ 
 satisfying conditions $(\alpha)_w-(\delta)_w$. Let $\ma_w$ be the 
subset of $R^o$  defined by the condition 
\begin{equation}
y_i 
\quad
{\rm divides}
\ \
 \Wr(\,y_i', \, T_i y_{i-1}y_{i+1})
\ \
{\rm for}\ \ i = 1, \dots ,  N\ ,
\end{equation}
Using the monicity of $y_i$ and long division, we can write the divisibility condition 
as a system of equations in the coefficients of $y_i$, $y_{i-1}$, and $y_{i+1}$. 
Hence $\mathcal{A}_w$ is a closed subscheme of $R_w^o$,

Consider the morphism 
$$
\Theta_w\ :\ \Fl_w^o\ \to\ R^o_w\ , 
\qquad
(V,E_{\bull})\ \mapsto\ (y_1, \dots , y_N)\
=\ (\,\Wr(\, E_1 ) / K_1, \dots , \Wr(\, E_N ) / K_N)\ .
$$ 
For $x\in \Fl^o_w$, condition $(\eta)_w$ holds and by Lemma ~\ref{new1}, 
$\Theta$ induces a morphism of schemes $\Theta:\Fl_w^o\to\ma_w$.

\begin{theorem}\label{mainte2}
The morphism \ $\Theta_w\,:\,\Fl^o_w\,\to\, \ma_w$\
is an isomorphism of schemes.
\end{theorem}

The proof of Theorem ~\ref{mainte2} is similar to that of Theorem
~\ref{mainte}. Notice that in Theorem ~\ref{Main}, no assumptions
were made on the degrees of $y_1, \dots , y_N$.

\bigskip

Consider the space $\tilde {R}_w = \prod_{i=1}^N \cee^{l^w_i}$ with
coordinates $(t^{(i)}_j)$, where $i=1,\dots,N$, \ $j=1,\dots,l^w_i$.
 The product of symmetric groups 
$\Sigma^w = \Sigma_{l^w_1}\times \dots\times \Sigma_{l^w_N}$
acts on $\tilde R$ by permuting coordinates with the same upper index.
Define a map
$$
\Gamma\ :\ \tilde {R}_w \ \to\ R\ , 
\qquad
(t^{(i)}_j)\ \mapsto \ (y_1,\dots,y_N)\ ,
$$  
where
$
y_i\ =\ \prod_{j=1}^{l_i^w}\ (x-t^{(i)}_j)\ .
$
Define the scheme $\tilde {\ma}_w$ by the condition $\tilde {\ma}_w =
 \Gamma^{-1}(\ma_w)$. The natural map $\tilde{\ma}_w\to\ma$ is finite and
 \'{e}tale.  The scheme $\tilde{\ma}_w$ is $\Sigma^w$-invariant. The
 scheme $\tilde{\ma}_w$ lies in the $\Sigma^w$-invariant subspace $\tilde
 {R}^o$ of all $(t^{(i)}_j)$ with the following properties for every $i$:
\begin{enumerate}
\item[$\bullet$] The numbers $t^{(i)}_1, \dots , t^{(i)}_{l^w_i}$ are distinct;
\item[$\bullet$]
The sets $\{t^{(i)}_1, \dots , t^{(i)}_{l^w_i}\}$
and $\{t^{(i+1)}_1, \dots , t^{(i+1)}_{l^w_{i+1}}\}$ do not intersect;
\item[$\bullet$]
The sets $\{t^{(i)}_1, \dots , t^{(i)}_{l^w_i}\}$
and $\{z_1, \dots , z_n\}$ do not intersect.
\end{enumerate}
The following lemma is proved in the same manner as Lemma ~\ref{oldone}.
\begin{lemma}${}$

\begin{enumerate}
\item[$\bullet$]
 The connected components of $\ma_w$ and $\tilde{\ma}_w$ are irreducible.
\item[$\bullet$]
  The reduced schemes underlying $\ma_w$ and $\tilde{\ma}_w$ are smooth. 
\item[$\bullet$]
 If $C$ is a connected component of $\ma_w$, then the group
$\Sigma^w$ acts transitively on the connected components of
$\Gamma^{-1}(C)$.
\end{enumerate}
\end{lemma}

Consider on $\tilde {R}^o_w$ the regular rational function 
$$
\Phi_w(t^{(i)}_j)\ =\ \prod_{i=1}^N \prod_{j=1}^{l^w_i} \ T_i(t^{(i)}_j)^{-1} \
\prod_{i=1}^{N-1} \prod_{j=1}^{l^w_i}\prod_{k=1}^{l^w_{i+1}}\
(t^{(i)}_j- t^{(i+1)}_k)^{-1}\
\prod_{i=1}^{N}\prod_{1\leq j< k \leq l^w_i}\
(t^{(i)}_j-t^{(i)}_k)^{2}\ .
$$ 
This $\Sigma^w$-invariant function is called {\it the master function
associated with the basic situation and the permutation $w\in \Sigma_{N+1}$}.

Define the scheme $\tilde {\ma}'_w$ as the subscheme in $\tilde R^o$
of critical points of the master function.
The following is a generalization of Lemma ~\ref{bethe eqn lemma} and
is proved in an identical fashion.
\begin{lemma}\label{ba11}
The subschemes $\tilde {\ma}_w$ and $\tilde {\ma}'_w$ of  $\tilde R^o$ coincide.
\end{lemma}
Let $(V,E_{\bull})\in \Fl^o_w$. Denote $\bs y =
\Theta(V,E_{\bull})\in\ma_w$. Pick a point $\bs t \in
\Gamma^{-1}(\bs y)$. Let $C$ be the unique irreducible component of
$\ma_w$ containing $\bs y$ and $\tilde{C}$ the unique irreducible
component of $\tilde{\ma}_w$ containing $\bs t$. The following is a
generalization of Theorem ~\ref{crit point multiplicity} with a
similar proof.
\begin{theorem}
The geometric multiplicity of the scheme $\Omega$ at $V$ equals
the geometric multiplicity of $\tilde{C}$.
\end{theorem}

\noindent
{\bf Example:} Let $N=1$, $n=3$, $z_1=1$,
$z_2=\omega$, $z_3=\omega^2$ where $\omega=e^{\frac{2\pi i}{3}}$.
Let $d=3$, $\aA(1)\,=\,\aA(\omega)\,=\,\aA(\omega^2)\,=\,\aA(\infty)\,=\,(1,0)$. 
Let $w$ be the transposition (12) in $\Sigma^2$.

It is easy to see that $T_1(x)=x^3-1$,\ {} $l^w_1=1$,   and
$\Phi(t)\,=\,{T_1(t)}^{-1}$. The critical scheme of 
 $\Phi$ is  $\{t\ |\ t^2\,=\,0\}$, namely $\spec(\cee[t]/(t^2))$. 
In other words, the master function
has one critical point at $t=0$
of multiplicity $2$. 

The polynomial $y_1$ associated to the critical point
is the polynomial $x$. The equation
$\Wr(\,y_1,\,\tilde{y_1})\,=\,T_1$ has solutions  $\tilde {y_1} 
\,=\, 1 + x^3/2 - cx$ with $c\in\cee$.
The associated two dimensional space of polynomials $V$\ is the $\C$-span of 
$x$ and $x^3 + 2$. The ramification points of $V$ are $1,\, \omega, \,\omega^2,\,
\infty$ with ramification sequences all equal to $(1, 0)$.

Counted with multiplicity there are two points of $\Omega$. The
associated cohomology product is the $4$th power of the hyperplane
class in $\Gr(2,4)$ which is $2$ points. It follows from Theorem
\ref{mainte2} that $\Omega$ is set-theoretically exactly one point $V$ counted 
with multiplicity 2. 

It is easy to see that $V$ admits a first order deformation in $\Omega$.
The deformation is given by
the $\C$-span of $x+\epsilon$ and $x^3 + 2 - 3 \epsilon\, x^2$
 over $\cee[\epsilon]/(\epsilon^2)$. Indeed we have
$$
\Wr\left(x+\epsilon,-1-\frac{x^3}{2}+\frac{3\epsilon\,x^2}{2}\right)\ =\
x^3 -1 \ {} \ {} \ {\rm mod}\, \epsilon^2\ .
$$

\section{Critical points of master functions}
Let $l_1,\dots,l_N$ be nonnegative integers and $z_1,\dots,z_n$
distinct complex numbers.  For $s=1, \dots , n$, fix nonnegative integers
$m_s(1), \dots , m_s(N)$. Define polynomials 
$T_1, \dots , T_N$ by the formula
$T_i\,  =\,  \prod_{s=1}^n \, (x-z_s)^{m_{s}(i)}$. 
{\it The master function $\Phi$ associated to this data} is the rational function 
$$
\Phi(t^{(i)}_j)\ =\ \prod_{i=1}^N \prod_{j=1}^{l_i} \ T_i(t^{(i)}_j)^{-1} \
\prod_{i=1}^{N-1} \prod_{j=1}^{l_i}\prod_{k=1}^{l_{i+1}}\
(t^{(i)}_j- t^{(i+1)}_k)^{-1}\
\prod_{i=1}^{N}\prod_{1\leq j< k \leq l_i}\
(t^{(i)}_j-t^{(i)}_k)^{2}
$$ 
of $l_1+\dots+l_N$ variables $t^{(1)}_1,\dots,t^{(1)}_{l_1},\dots,
t^{(N)}_{1},\dots,t^{(N)}_{l_N}$ considered on the set of points where
\begin{enumerate}
\item[$\bullet$] The numbers $t^{(i)}_1, \dots , t^{(i)}_{l_i}$ are distinct;
\item[$\bullet$]
The sets $\{t^{(i)}_1, \dots , t^{(i)}_{l_i}\}$
and $\{t^{(i+1)}_1, \dots , t^{(i+1)}_{l_{i+1}}\}$ do not intersect;
\item[$\bullet$]
The sets $\{t^{(i)}_1, \dots , t^{(i)}_{l_i}\}$
and $\{z_1, \dots , z_n\}$ do not intersect.
\end{enumerate}
Let $l_0=l_{N+1}=0$. Define 
$c_i\ =\ i-1\ +\ l_i - l_{i-1}\ +\ \sum_{j=1}^{i-1}\sum_{s=1}^n\ m_{s}(j)$\
for \  $i = 1 ,\dots , N+1.$
\begin{lemma} If either  $c_i<0$ for some $i\, =\, 1,\dots,N+1$ \ or\ 
$c_i\, =\, c_j$ for some  $1\, \leq i\, <\,  j\leq\,  N+1$, then
 $\Phi$ has no critical points.
\end{lemma}
\begin{proof}
Let $(t^{(i)}_j)$ be  a critical point of
$\Phi$. Set $y_i\ =\ \prod_{j=1}^{l_i}\ (x-t^{(i)}_j)$ for
$i=1,\dots,N$. Then $y_i$ divides $\Wr(\,y_i',\,T_i\,
y_{i-1}\,y_{i+1})$, where as usual we set $y_{N+1}=y_0=1$. 

We have $c_i\ =\ i-1\ + \ \deg\,y_i\ -\ \deg\,y_{i-1}\ +\
\sum_{j=1}^{i-1}\,\deg\,T_j$.  By part (3) of Theorem ~\ref{Main},
$c_1, \dots , c_{N+1}$ are pairwise distinct non-negative
numbers.
\end{proof}

Assume  $c_1, \dots , c_{N+1}$ are pairwise distinct nonnegative integers. 
 Let $u \in \Sigma_{N+1}$ be the
permutation such that $c_{u(1)}\, >\, c_{u(2)}\, >\, \dots\, >\,
c_{u(N+1)}$.  Let $d = \max\,\{c_1,\dots,c_{N+1}\}$.

For $z\in
\{z_1,\dots,z_n\}$ define the ramification sequence $\aA(z)$ by the rule
$$
a_{j}(z)\ =\ \sum_{\ell=j}^{N}\, m_{N+1-\ell}(z)\ ,
\qquad
 j = 1 , \dots , N+1\ .
$$
 Define the ramification sequence $\aA(\infty)$ at
$\infty$ by the rule
$$
a_{j}(\infty)\ =\ d\ -\ (N+1)+j\ -\ c_{u(j)}\ , 
\qquad j = 1 , \dots , N+1\ .
$$ 
\begin{proposition}  
The master function $\Phi$ of this section
is the same as the master function
of the basic situation of Section ~\ref{extendedbasic} 
associated with the space
$W = \C_d[x]$,\ ramification sequences
$\aA(z_1),\,\dots,\,\aA(z_n),\,\aA(\infty)$,\ and the permutation
$u^{-1}$ in $\Sigma_{N+1}$. 
\hfill
$\square$
\end{proposition}

\appendix\section{}

\subsection{Equation $\Wr(\,y ,\, \tilde{y})\,=\,T$}
\begin{lemma}\label{BA}
Let $A$ be an algebra over $\cee$. Let $y, T\in A[x]$ . 
\begin{enumerate}
\item If the equation $\Wr(\,y , \tilde{y})\,=\,T$ has a solution
$\tilde{y}\in A[x]$ then $y$ divides $\Wr(\,y',T)$.
\item If the ideal $(y,y')\ =\ A[x]$ and $y$ divides $\Wr(\,y',T)$, then the equation $\Wr(\,y , \tilde{y})\,=\,T$ has a solution
$\tilde{y}\in A[x]$.
\end{enumerate} 
\end{lemma}

\begin{proof} 
If $T\,=\,\Wr( \,y , \tilde{y})$, then
$\Wr(\,y' , \Wr(\,y , \tilde{y}))\,=\,
y''\,(y'\tilde{y}\,-\,y\,
\tilde{y}')\,-\,y'\,(y''\,\tilde{y}\,-\,y\,\tilde{y}'')\ =\ 
y\,(y'\,
\tilde{y}''\,-\,y''\,\tilde{y}')$, and $y$ divides
$\Wr(\,y' , T)$.

  For the other direction, write
$T = a y + b y'$ for suitable polynomials $a, b$. Then
\linebreak
$T = ( a + b')\,y \,+\, \Wr(\,y , b)$. 
We have $\Wr(\,y', T) = \Wr(\,y', y(a + b')) + \Wr(\,y', \Wr(\,y, b))$.
The last term is divisible by $y$. If $y$ divides $\Wr(\,y', T)$, then
$y$ divides 
$ \Wr(\,y' ,\,y\,(a + b'))$. 
This implies that $y$ divides
$(y')^2(a + b')$. 
Writing $c\,y'\,=\,1\,-\,d y$ for suitable polynomials $c$ and $ d$,
we deduce that $y$ divides $(a + b')$. Thus, $y(a + b') = y^2 e$
for a suitable polynomial $e$. Let $f$ be a
polynomial whose derivative is $-e$. Then $\Wr(\,1 , f)\ =\ e$, and
$\Wr(\,y , y f)\ =\  y(a+b')$. Finally,
$\Wr(\,y , y f+ b) = y(a+b') + \Wr(\,y, b) = T$.
\end{proof}

\begin{lemma}\label{divisionremark}
 Suppose $y, f \in A[x]$ and $y$ is monic. 
\begin{enumerate}
\item There exist unique  $q, r \in A[x]$ such that
$f\,=\,y\,q\,+\,r$ and $\deg \,r\, <\deg\,y$. 
\item Suppose that $y\,=\,\prod_{j=1}^m\, (x-t_j), t_j\in A$ and $t_j-t_k$ are units in $A$
for $j\neq k$.  Then the ideal in $A$ generated by coefficients of
the polynomial $r$ in (1), coincides with the ideal in $A$ generated by
$f(t_1), \dots , f(t_m)$.
\end{enumerate}
\end{lemma}
\begin{proof}
Existence and uniqueness of $q$ and $r$ follows from the monicity of
$y$ and long division. 

For the second part, write
\ $f(t_j)\ =\ r_{m-1}t_j^{m-1}\ + \ \dots \ +\ r_0$.\ 
From Cramer's rule and the formula for the Vandermonde determinant we deduce that
the coefficients of $r$ lie in the ideal generated by $f(t_1), \dots , f(t_m)$.
%
\end{proof}

\begin{lemma}\label{mon1}
Let  $(A,\emm)$ be a local ring and $y_1, y_2\in A[x]$
with $y_1$ monic. Suppose that the ideal generated by the  reductions 
$(\bar{y}_1, \bar{y}_2)=(A/\emm)[x]$. Then the ideal $(y_1,\,y_2)$ equals $A[x]$.
\end{lemma}

\begin{proof}
The $A$-module $A[x]/(y_1)$ is a finite $A$-module because
$y_1$ is monic. The quotient $M\,=\,A[x]/(y_1,y_2)$
 of $A[x]$ by the ideal $(y_1,y_2)$ is therefore a finite $A$-module.  Now, $M\tensor
{A/\emm}\, =\, ({A/\emm})[x]\,/\,(\bar{y}_1,\bar{y}_2)\,=\,0$ by hypothesis. 
By Nakayama's Lemma \,(~\cite{mat}, Theorem 4.8), we 
conclude that $M\,=\,0$.
\end{proof}

\begin{lemma}\label{equationes} Let $A$ be an algebra, 
$y  = \prod_{j=1}^m\,(x-t_j)\in A[x]$ 
with $t_j\in A$ and $T = \prod_{l=1}^n\,(x-z_l)\in A[x]$ with $z_l\in A$. Then
\begin{enumerate}
\item The ideal $(y,\,y')=A[x]$ if and only if\ $t_i-t_j$ are units  in $A$ for all
$i< j $. 
\item The ideal $(T,\,y)=A[x]$ if
and only if \ $t_j-z_l$ are units in $A$ for all $j, l$.
\item Assume that the ideals $(y,\,y')= A[x]$ and $(T,\,y)=A[x]$. Then,
$y$ divides
$\Wr(\,y', T)$ if and only if the following system of equations holds:
\begin{equation}\label{crit}
\sum_{l\in\{1,\dots,m\}-j}\,
\frac{2}{t_j-t_l}\ -\ 
\sum_{l=1}^n\,\frac{1}{t_j-z_l}\ =\ 0\ ,
\qquad
 j = 1, \dots , m\ .
\end{equation}
\end{enumerate}
\end{lemma}
Notice that the system of equations (~\ref{crit}) coincides with the  critical point equations for the function
$$
\Phi(t_1,\dots,t_m)\ = \!\!
\prod_{1\leq i<j\leq m}\,(t_i-t_j)^2\
\prod_{j=1}^m\prod_{l=1}^n \,(t_j-z_l)^{-1}
\ =\!\! \prod_{1\leq i<j\leq m}\,
(t_i-t_j)^2\ 
\prod_{j=1}^n  T(t_j)^{-1}\ 
$$
of variables $t_1, \dots , t_m$.

\begin{proof} 
If the ideal $(y,y')\,=\,A[x]$, write $1 = ay+by'$ and substitute $x = t_i$ to conclude
that $t_i - t_j$ is invertible for $i \neq j$. 

For the other direction, let
$M=A[x]/(y,y')$. Clearly, $M/\emm M = ({A/\emm})[x]/(\bar{y},\bar{y'}) = 0$,\ 
which is guaranteed by the assumption and the standard theory of fields. 
Since $y$ is monic, by Lemma ~\ref{mon1}, $M=0$. The assertion (2) is proved similarly.

It is easy to see that (3) follows from Lemma ~\ref{divisionremark}
 and Lemma ~\ref{BA}.
\end{proof}

\begin{lemma}\label{mon2}
For $y_1,y_2, T \in A[x]$, assume that  the ideal $(y_1,\,y_2)\ =\ A[x]$
and $y_1$ divides $y_2 \,T$. Then $y_1$ divides $T$.
\end{lemma}
\begin{proof}
Write $1 = a y_1 + b y_2$. Hence $T = a T y_1 + b T y_2$. 
The last two terms are divisible by $y_1$. Hence
 $y_1$ divides $T$.
\end{proof}

\subsection{Proof of Lemma ~\ref{setup}}\label{setupproof}
The implication $(1)\Rightarrow (3)$ is immediate.

If (2) holds , then $V$ is finitely generated module which is a direct
summand of a free module. This implies that it is free (being a direct
summand implies that it is projective and projective modules over
local rings are free). The morphism $V\tensor A/\emm\to \cee[x]$ is a
direct summand of the isomorphism $A[x]\tensor A/\emm\to \cee[x]$ and 
hence is injective. This gives (1).

We prove $(1)$ and $(2)$ assuming $(3)$.  Pick a large
integer $d$ so that $V\subset A_{ d}[x]$.  By assumption we can
find $v_1,\dots,v_{d+1-k}\in A_{ d}[x]$ such that the collection
$\bar{u}_1, \dots , \bar{u}_k,$ \linebreak $ \bar{v}_1,$
$\dots , $ $\bar{v}_{d+1-k}$ form a
basis for the $\cee$-vector space $\cee_{ d}[x]$. This implies
that the change of basis matrix from the standard basis of $A_{
d}[x]$ to $u_1,\dots, u_k,v_1,\dots,v_{d+1-k}$ has the determinant that
does not vanish upon reduction to the residue field. 
Therefore, the  determinant is  a unit in $A$. Hence
$u_1,\dots,u_k,v_1,\dots, v_{d+1-k}$ form a free basis in $A_{
d}[x]$. This proves $(1)$ and $(2)$, where for (2) we let
$M \,=\, \spanA(v_1,\dots,v_{d+1-k})\oplus x^{d+1}A[x]\subset A[x]$.
\subsection{Multiplicity}
We will recall the algebro-geometric definitions of multiplicity from
\linebreak ~\cite{int}, Section 1.5. In this section we will need to
consider local rings whose residue field may be  different from $\cee$.

Let $X$ be an irreducible algebraic scheme. The geometric multiplicity
of $X$, denoted by $\mm(X)$, is the length of the local ring of $X$ at its
generic point. Explicitly, if $\spec(A)$ is an affine open subset of
$X$, then $A$ has exactly one minimal prime ideal. Denote it by
$\pp$. The localisation $A_{\pp}$ is an Artin ring and is therefore of
finite length. The integer $\mm(X)$ is the length of $A_{\pp}$. A more
practical way of computing $\mm(X)$ is obtained from Proposition
~\ref{transversalmultiplicity}.

{\bf Example:} Consider the geometric multiplicity of
the ``doubled line'' \linebreak
$\spec(\cee[x,y]/(x^2))$. This has exactly one
minimal prime ideal, namely $(x)$. The localisation at this minimal prime
ideal is the ring $\cee(y)[x]/(x^2)$, which is of length $2$.

Now, we will discuss properties of geometric multiplicity, 
linking it to ``multiplicity in the transversal direction''. The first
property is (see ~\cite{int}, Example A.1.1):

\begin{proposition} 
If $X$ is an irreducible $0$-dimensional scheme (i.e a fat point),
then $\mm(X)$ is the dimension over $\Bbb{C}$ of the ring of functions
of $X$,
$$
\mm(X)\ =\
\dim_{\Bbb{C}}\,\Gamma(X,\mathcal{O}_X)\ .
$$
\end{proposition}
In general, if $X$ is an irreducible scheme, then $X$ is reduced if
and only if its geometric multiplicity is $1$.

The following proposition follows from Lemma 1.7.2 in ~\cite{int}.
\begin{proposition} \label{transversalmultiplicity}
Suppose that $X$ is an irreducible subscheme of 
$\Bbb{C}^n$. Let $X_{red}$ be the reduced subscheme corresponding
to $X$. The subscheme
$X_{red}$ can be considered to be a closed subscheme of
$X$. Let $U$ be the smooth locus of $X_{red}$. Let $H$ be a hyperplane
in $\Bbb{C}^n$ which meets $U$ transversally at a point $x\in U$. Let
$D$ be the irreducible component of $X\cap H$ which contains $x$
(there is exactly one such irreducible component). Then,
$$
\mm(X)\ =\ \mm(D)\ .
$$
Iterating this procedure, we obtain the following statement. 
Suppose $T$ is a plane in
$\Bbb{C}^n$ of dimension complementary to $\dim\,X$, which meets $U$
transversally at a point $x\in U$ (there could be other points of
intersection). Then, the multiplicity of $X$ equals to the dimension
over $\Bbb{C}$ of the localization at $x$ of the algebra of functions
on the scheme $X\cap H$.
\end{proposition}

There is one other standard property of multiplicity that we will need. 
Recall that a smooth morphism between schemes is a flat morphism with smooth fibers. 
\begin{proposition}\label{pga} 
Let $f:X\to Y$ be a smooth morphism 
between irreducible schemes. Then, $\mm(X)=\mm(Y)$.
\end{proposition}
\begin{proof}
We will use the notations and definitions of ~\cite{int}.
  Let $X_{\red}$ and $Y_{\red}$ be the reduced schemes
underlying $X$ and $Y$.  Then by definition $[X]=\mm(X)[X_{\red}]$ and
$[Y]=\mm(Y)[Y_{\red}]$.  The smoothness of $f$ tells us that
$f^{-1}(Y_{\red})$ is reduced and hence
$f^{-1}(Y_{\red})=X_{\red}$. Clearly $f^{-1}(Y)=X$. Apply Lemma 1.7.1
in ~\cite{int} to see that \linebreak
$f^*[Y]=[f^{-1}(Y)]$.  This gives
$\mm(Y)[X_{\red}] =  \mm(X)[X_{\red}]$ and therefore $\mm(X) = \mm(Y)$.
\end{proof}
\subsection{Multiplicity in intersection theory}
Irreducible subvarieties $X_1,\dots,X_r$ of a smooth variety $X$ are
said to intersect properly, provided each irreducible component of
$X_1\cap\dots\cap X_r$ is of dimension
$\dim\,X\,-\,\sum_{j=1}^r\,(\dim\,X\,-\,\dim\,X_j\,)$. We will use the following
basic result.

 Denote the smooth locus of $X_i$ by $X_i^o$. Suppose that
$X_1, \dots , X_r$ intersect properly in a finite set,\ that is, the
expected dimension of the intersection is $0$. Suppose that
 $$
X_1\ \cap\ \dots\ \cap \ X_r\ =\ X_1^o\ \cap\ \dots\ \cap X_r^o \ .
$$ 
Then we  have an equality of cohomology classes
in $ H^*(X)$, 
$$
\prod_{i=1}^r\,[X_i]\ =\ c\ [\pt]\ ,
$$ 
where $c$ is the sum of the
multiplicities of the irreducible components 
of the scheme theoretic intersection 
$X_1\cap\dots\cap X_r$ and $[\pt]$ the class of a point. 

This statement follows from ~\cite{int}, Proposition 7.1. 

\subsection{Standard results in the theory of schemes}
For a scheme $X$ and a $\cee$-algebra $A$, we let 
$$
X(A)\ =\ \home( \spec(A) , X)\ .
$$ If $A$ is a local ring, and $s\in X(A)$, then we denote the induced
point in $X(\cee) = X({A/\emm})$ by $\bar{s}$. If $x\in X({A/\emm})$ is
given, we let $X_x(A) = \{s\in X(A)\ |\ \bar{s}=x\}$.  For $s\in
X_x(A)$ there corresponds a local homomorphism of local rings
$\mathcal{O}_{X,x}\to A$, where $\mathcal{O}_{X,x}$ is the local ring
of $X$ at $x$.

\begin{lemma}\label{schema}
Let $f:X\to Y$ be a finite morphism of schemes. Then, $f$ is a closed
immersion if and only if for every local ring $A$, the induced mapping
$X(A)\to Y(A)$ is injective.
\end{lemma}
\begin{proof}
If $f:X\to Y$ makes $X$ a subscheme of $Y$, then clearly
 $X(A)\to Y(A)$ is injective.

To go the other way, let $x\in X$ and $y=f(x)$. By taking $A=\cee$, we
see that $f^{-1}(y)$ is the singleton $\{x\}$. Denote the local ring
of $X$ at $x$ by $(\mathcal{O}_{X,x},m_x)$ and that of $Y$ at $y$ by
$(\mathcal{O}_{Y,y},m_y)$.

Now let $A=\cee[\epsilon]/(\epsilon^2)$. Consider the induced mapping
$X_x(A)\to Y_y(A)$. This is once again injective by hypothesis. It is
a basic fact that $X_x(A)=\home(m_x/m_x^2,\cee)$ and
$Y_y(A)=\home(m_y/m_y^2,\cee)$. Therefore the hypothesis implies that
$\home(m_x/m_x^2,\cee)\to\home(m_y/m_y^2,\cee)$ is injective or that
the natural morphism $m_y/m_y^2\to m_x/m_x^2$ is surjective.  By
~\cite{hartshorne}, II.7.4, we conclude that the map
$\mathcal{O}_{Y,y}\to \mathcal{O}_{X,x}$ is surjective. This shows
that $X\to Y$ is a closed inclusion.
\end{proof}

\begin{lemma}\label{schema2}
Let $f:X\to Y$ be a closed immersion of schemes. Then $f$ is an
isomorphism if and only if $f:X(A)\to Y(A)$ is surjective for every
local ring $A$.
\end{lemma}
\begin{proof}
If $f$ is an isomorphism, then $X(A)\to Y(A)$ is clearly bijective for
any local ring $A$.

To go the other way, let $x\in X$ and $y=f(x)$. Denote the local ring
of $X$ at $x$ by $(\mathcal{O}_{X,x},m_x)$ and that of $Y$ at $y$ by
$(\mathcal{O}_{Y,y},m_y)$.

Suppose that $y$ is in an open affine subset $\spec(B)\subseteq
Y$. Let $\emm$ be the maximal ideal of $\spec(B)$ corresponding to the
point $y$. Then $\mathcal{O}_{Y,y}$ is the same as the localization
$B_m$ of $B$ at $m$. The natural map $B\to B_m$ gives a point $\eta$
in $\spec(B)(\mathcal{O}_{Y,y})\subseteq
Y(\mathcal{O}_{Y,y})$. Clearly $\bar{\eta}=y$ and therefore $\eta\in
Y_y(\mathcal{O}_{Y,y})$. The map of local rings corresponding to
$\eta$ is the identity map $\mathcal{O}_{Y,y}\to \mathcal{O}_{Y,y}$

By the given hypothesis, there exists $\theta\in X(\mathcal{O}_{Y,y})$
such that $f(\theta)=\eta$. The reduction of this point is $x\in
X(\cee)$ because the reduction has to sit over $y\in Y(\cee)$.

Therefore we obtain a diagram
$$
\xymatrix{\mathcal{O}_{Y,y}\ar[r]^{f}\ar[rd]^{=} & \mathcal{O}_{X,x}\ar[d]^{\theta}\\
                                         & \mathcal{O}_{Y,y}}
$$

Hence $\mathcal{O}_{Y,y}\to \mathcal{O}_{X,x}$ is both injective and
surjective whence an isomorphism.
\end{proof}

The following proposition is standard and follows from the description of
 Schubert varieties as degeneracy loci (cf. ~\cite{int}, chapter 14).
\begin{proposition}\label{functor}
Let $A$ be a local ring and $W$ a $\cee$-vector space of rank
$d+1$. Then, $\Gr(N+1,W)(A)$ is the set of free submodules $V\subset
W\tensor A$ of rank $N+1$, such that $V\tensor A/\emm \to W$ is
injective.

The subset $\Omega^o_{\aA}(\mf)(A)\subseteq \Gr(N+1,W)(A)$ consists of
 submodules $V$ such that there exists an $A$-basis 
$u_1,\dots,u_{d+1}$ of $W\tensor A$
with the following properties:
\begin{enumerate}
\item[$\bullet$] $F_i$ is the $A$-span of $u_1, \dots , u_i$ for $i  =  1,
\dots , d+1$;
\item [$\bullet$]
$V$ is the $A$-span of the elements $u_{d-N+j-a_j}$ for $j=1,\dots,N+1$.
\end{enumerate}
\end{proposition}

\bibliographystyle{plain}
\def\noopsort#1{}

\end{document}